\documentclass[a4paper,11.5pt,reqno]{amsart}

\usepackage[pdftex]{graphicx}
\usepackage{amsmath,amssymb,amsthm}
\usepackage{enumitem}
\usepackage{xcolor}
\usepackage{hyperref}
\usepackage{lineno}

\numberwithin{equation}{section}

\hypersetup{colorlinks = true, linkcolor = blue, urlcolor = blue, citecolor = blue}

\usepackage{lineno}
\modulolinenumbers[1]

\theoremstyle{plain}
\newtheorem{theorem}{Theorem}[section]
\newtheorem{conjecture}[theorem]{Conjecture}

\newtheorem{proposition}[theorem]{Proposition}

\theoremstyle{definition}

\theoremstyle{remark}

\begin{document}

\title[Cremona Maps and Krall-Jacobi Polynomials]
      {Two Families of Cremona Maps and orthogonal Krall-Jacobi Polynomials}

\author[H. Ruhland]{Helmut Ruhland}
\address{Santa F\'{e}, La Habana, Cuba}
\email{helmut.ruhland50@web.de}



\begin{abstract}
Two infinite families of Cremona maps depending on one real parameter are given. For all integers
$n \ge 1$ the first family of Cremona maps consists of group elements in $Bir \left( \mathbb{P}^{n} \right)$
with bidegree $(n, n)$, the second family of Cremona maps consists of group elements in $Bir \left( \mathbb{P}^{2n} \right)$ with bidegree $(2 n, 2 n)$. For the first family and $n \ge 5$, for the second
family and $n \ge 3$ the existence of this group elements and the properties depend on a conjecture.
But computational results suggest that the conjecture is true for all $n$.
\end{abstract}

\date{\today}

\maketitle

\section{Introduction}

Among the known infinite families of Cremona maps is the family of the standard involutions
$\sigma_n : \mathbb{RP}^{n} \rightarrow \mathbb{RP}^{n}$ in $Bir \left( \mathbb{P}_\mathbb{R}^{n} \right)$
for all ranks of the Cremona groups. This family has a very simple definition, see e.g. \cite{Des2012}. \\

A goal of this article is to present two further, infinite families of Cremona maps, maybe novel,
depending on \emph{one real parameter} for all integers $n \ge 1$:
\begin{itemize}
   \item the first family of Cremona maps in Cremona groups of all ranks consists of group elements in
         $Bir \left( \mathbb{P}^{n} \right)$ with bidegree $(n, n)$.
   \item the second family of Cremona maps in Cremona groups only for even rank consists of group
         elements in $Bir \left( \mathbb{P}^{2n} \right)$ with bidegree $(2 n, 2 n)$.
\end{itemize}
The existence and the properties of these maps depend on two conjectures. I don't know a rigorous mathematical proof of the conjectures, even not in the $6$ cases for the maps in the Cremona groups
$Bir \left( \mathbb{P}^{4} \right)$ and lower rank. Instead of presenting proofs, a second
goal of this article is to show \emph{strong computational evidence} for the conjectures to be true
by symbolic computations with a computer algebra system. \\

The construction of these two families is based on semi-classical orthogonal polynomials that were already
introduced by Krall (1940) \cite{Krall1940}, \cite{Koor1984} and \cite{ArMaAl2002}.
In newer publications, see e.g. \cite{Dur2022}, see formulae (1.4) ... (1.6) these polynomials belong
to the so named orthogonal \emph{Krall-Jacobi families}. \\
The construction of the maps uses the fact, that the orthogonal Krall-Jacobi polynomials used in
\emph{this article} can be expressed as a sum of classical, symmetric Jacobi polynomials. \\

If the reader is only interested in the definition of the maps and in the conjectures and not
in computational details all appendices can be ignored. \\

These 2 families of Cremona maps play a role in Numerical Analysis for quadrature rules of splines,
see appendix \ref{SectNumAna}.

\section{The construction of the first family of 1-parameter Cremona maps}

Define a vector $\textbf{l} = (l_0, l_1, \dots , l_{n-2}, l_{n-1}) \in \mathbb{R}^n$.
Orthogonal polynomials can be defined by their weight function (here a distribution):
\begin{equation} W_Q(\textbf{l}, x) = (1 - x)^{n} \left( 1 + \sum_{i=0}^{n-1} l_i \delta^{(i)}(1+x) \right) \label{WeightQ} \end{equation}
The $\delta^{(i)}()$ are a Dirac distribution and its $i^{th}$ derivations. Because the Dirac's in this
weight are located on the left side of the interval $[-1, +1]$ I named the vector $\textbf{l}$ like (l)eft
and the corresponding orthogonal polynomials one-sided. \\

With the so defined weight function the \textit{semi-classical} one-sided orthogonal Jacobi polynomial
$Q_k  (\textbf{l}, x)$ of degree $k$ is defined by this property:
\begin{align}
\enspace & \int\limits_{-1}^{+1} W_Q(\textbf{l}, x) \; Q_k (\textbf{l}, x) \; Q_l (\textbf{l}, x) \, dx = 0 \quad  \hbox{for} \enspace k \ne l  \label{OrthQ}
\end{align}
$Q_k (\textbf{0}, x) = P_k^{(n,0)}(x)$ is a classical orthogonal Jacobi polynomial. \\

These \textit{semi-classical} Jacobi polynomials $Q_k (\textbf{l}, x)$ can be expressed by a sum of
\textit{symmetric} Jacobi polynomials $J_{k} (x) = P_k^{(n,n)} (x)$, the coefficients $j_{k,i} (\textbf{l})$ are polynomials in the components of the vector $\textbf{l}$ and $k$:
\begin{equation} Q_k (\textbf{l}, x) = \sum_{i=0}^{n} j_{k,i} (\textbf{l}) J_{k-i} (x) \label{QasSumJ} \end{equation}

The $j_{k,i} (\textbf{l})$ in (\ref{QasSumJ}) now define a polynomial map:
\begin{align}
f_k : \enspace  & \mathbb{R}^n & \rightarrow & \enspace \mathbb{RP}^{n} \nonumber \\
& \textbf{l} = ( l_0, \, \dots , l_{n-1} ) & \mapsto & \enspace ( j_{k,0} (\textbf{l}) : \dots : j_{k,n} (\textbf{l}) ) \label{MapTrans}
\end{align}
By definition $k \in \mathbb{N}_+$, we can extend the domain to $k \in \mathbb{R}$. The components of the
maps are then polynomial functions in $k$. To distinguish the map $f_k$ only defined for integer indices $k$
from the map defined for $k \in \mathbb{R}$ we write the later maps with an argument instead of an index
as $f (k)$.  

\begin{conjecture} \label{ConjOneSided}
The map $f^* (k) : \mathbb{RP}^{n} \rightarrow \mathbb{RP}^{n}$, the homogenized $f (k)$ is birational.
This map and its inverse $f^{*-1} (k)$ have bidegree $(n, n)$. So it is a group element in the
Cremona group $Bir \left( \mathbb{P}_\mathbb{R}^{n} \right)$, the group of birational maps in the projective space $\mathbb{RP}^{n}$. The determinant of the Jacobian:
\begin{equation}
  \mathrm{det}\,\mathrm{jac} \, (f^* (k)) = D_{2n-1} (k) \quad P^{n-1} \; S_{n-1}^n
  \quad \text{of degree } n^2 - 1
\label{fForm} \end{equation}
is the union of \textbf{two components}: a hyperplane $P$ of multiplicity $n - 1$ and a hypersurface
$S_{n-1}$ of degree $n - 1$ and multiplicity $n$.
Let the homogenizing variable be $h$. Then the hyperplane has this form $P : h = 0$
and the hypersurface depends on $h$ and $l_1 \dots l_{n-1}$ but not on $l_0$. \\
$D_{2n-1} (k)$ is a polynomial of degree $2n - 1$ in $k$ and determines where 
the map $f^* (k)$ has codimension $1$ ($\mathrm{det}\,\mathrm{jac} = 0$). All zeros of this polynomial
occur for $k \in D = \{ -n \enspace (1/2) \enspace -1 \}$, in this set 2 consecutive elements
have an increment of $1/2$. So in this case the map is not birational.
For $n \ge 3$ the hypersurface $S_{n-1}$ is reducible for $k = -2n + 1, \, \dots, n - 2$.
The degree in $k$ of this map and its inverse is $2n^2$.
\end{conjecture}
Here the conjecture for the inverse maps:
\begin{conjecture} \label{ConjOneSidedInv}
The determinant of the Jacobian $\mathrm{det}\,\mathrm{jac}$ of $f^{*-1} (k)$
\begin{equation}
  \mathrm{det}\,\mathrm{jac} \, (f^{*-1} (k)) \sim P^{n(n-1)} \; S_{n-1}
  \quad \text{of degree } n^2 - 1
\label{fForm} \end{equation}
is the union of \textbf{two components}: a hyperplane $P$ of multiplicity $n (n - 1)$ and a hypersurface
$S_{n-1}$ of degree $n - 1$ and multiplicity $1$.
\end{conjecture}
For $n = 1, 2, 3$ and $4$ the 2 conjectures above are true, see appendices \ref{SectBir1}, \ref{SectBir2},
\ref{SectBir3} and \ref{SectBir4}. About the computation of these maps see appendix \ref{SectCAS}. \\
Further evidence for this conjecture is based on observations doing symbolic calculations with a CAS
for the case $n = 5$ and \emph{some low} $k$.

\section{The construction of the second family of 1-parameter Cremona maps}

Define two vectors $\textbf{l} = (l_0, l_1, \dots , l_{n-2}, l_{n-1}) \in \mathbb{R}^n$ and
$\textbf{r} = (r_0, r_1, \dots , r_{n-2}, r_{n-1}) $ $\in \mathbb{R}^n$.
Orthogonal polynomials can be defined by their weight function (here a distribution):
\begin{equation} W_M(\textbf{l}, \textbf{r}, x) = 1 + \sum_{i=0}^{n-1} l_i \delta^{(i)}(1+x) + \sum_{i=0}^{n-1} r_i \delta^{(i)}(1-x) \label{WeightM} \end{equation}
The $\delta^{(i)}()$ are a Dirac distribution and its $i^{th}$ derivations. Because the Dirac's in this
weight are located on the left and right side of the interval $[-1, +1]$ I named the vectors
$\textbf{l}, \textbf{r}$ like (l)eft and (r)ight and the corresponding orthogonal polynomials two-sided. \\
This weight function has the following involution as symmetry, it is called here left-right symmetry:
\begin{equation}
LR : x \mapsto -x \quad l_i \mapsto r_i \quad r_i \mapsto l_i  \quad
 \hbox{for} \enspace 0 \le i \le n - 1   
\label{LR_Sym} \end{equation}

With the so defined weight function the \textit{semi-classical} orthogonal two-sided Jacobi polynomial
$M_k  (\textbf{l}, \textbf{r}, x)$ of degree $k$ is defined by this property:
\begin{align}
\enspace & \int\limits_{-1}^{+1} W_M (\textbf{l}, \textbf{r}, x) \; M_k (\textbf{l}, \textbf{r}, x) \;
 M_l (\textbf{l}, \textbf{r}, x) \, dx = 0  \quad  \hbox{for} \enspace k \ne l \nonumber
\end{align}
$M_k (\textbf{0}, \textbf{0}, x) = P_k(x)$ is a classical orthogonal Legendre polynomial. The
left-right symmetry (\ref{LR_Sym}) results in the identity $M_k (\textbf{l}, \textbf{r}, x) =
 M_k (\textbf{r}, \textbf{l}, -x)$. \\

These \textit{semi-classical} Jacobi polynomials $M_k (\textbf{l}, \textbf{r}, x)$ can be expressed by
a sum of \textit{symmetric} Jacobi polynomials $J_{k} (x) = P_k^{(n,n)} (x)$, the coefficients $j_{k,i} (\textbf{l}, \textbf{r})$ are polynomials in the components of the vectors $\textbf{l}, \textbf{r}$ and $k$:
\begin{equation} M_k (\textbf{l}, \textbf{r}, x) = \sum_{i=0}^{2n} j_{k,i} (\textbf{l}, \textbf{r}) J_{k-i} (x) \label{MasSumJ} \end{equation}

The $j_{k,i} (\textbf{l}, \textbf{r})$ in (\ref{MasSumJ}) now define a polynomial map:
\begin{align}
g_k : \enspace  & \mathbb{R}^{2n} & \rightarrow & \enspace \mathbb{RP}^{2n} \nonumber \\
& \textbf{l}, \textbf{r} = ( l_0, \, \dots , l_{n-1}, \; r_0, \, \dots , r_{n-1} ) & \mapsto & \enspace ( j_{k,0} (\textbf{l}, \textbf{r}) : \dots : j_{k,2n} (\textbf{l}, \textbf{r}) ) \label{MapTrans}
\end{align}
The left-right symmetry (\ref{LR_Sym}) in the weight function acts on the $\textbf{j}$ components as:
\begin{equation}
LR : j_i \mapsto (-1)^i \, j_i \quad \hbox{for} \enspace 0 \le i \le 2n   
\label{LR_Sym_j} \end{equation}
By definition $k \in \mathbb{N}_+$, we can extend the domain to $k \in \mathbb{R}$. The components of the
maps are then polynomial functions in $k$. To distinguish the map $g_k$ only defined for integer indices $k$
from the map defined for $k \in \mathbb{R}$ we write the later maps with an argument instead of an index
as $g (k)$.  

\begin{conjecture} \label{ConjTwoSided}
The map $g^* (k) : \mathbb{RP}^{2n} \rightarrow \mathbb{RP}^{2n}$, the homogenized $g (k)$ is birational.
This map and its inverse $g^{*-1} (k)$ have bidegree $(2 n, 2 n)$. So it is a group element in the
Cremona group $Bir \left( \mathbb{P}_\mathbb{R}^{2n} \right)$. The determinant of the Jacobian:
\begin{equation}
  \mathrm{det}\,\mathrm{jac} \, (g_k^*) = D_{2n-1} (k) \quad P^{2n-1} \; S_{2n-1}^n \; \overline{S}_{2n-1}^n
  \quad \text{of degree } (2n)^2 - 1
\label{fForm} \end{equation}
is the union of \textbf{three components}: a hyperplane $P$ of multiplicity $2n - 1$ and two hypersurfaces
$S_{2n-1}$ and $\overline{S}_{2n-1}$ of degree $2n - 1$ and each of multiplicity $n$. These hypersurfaces
$S_{2n-1}$ and $\overline{S}_{2n-1}$ are permuted by the left-right symmetry (\ref{LR_Sym}). \\
$D_{2n-1} (k)$ is a polynomial of degree $2n - 1$ in $k$ and determines where 
the map $g^* (k)$ has codimension $1$ ($\mathrm{det}\,\mathrm{jac} = 0$). All zeros of this polynomial
occur for $k \in D = \{ -(2n-1)/2 \enspace (1) \enspace (2n-3)/2 \}$, in this set 2 consecutive
elements have an increment of $1$. So in this case the map is not birational.
For $n \ge 2$ the hypersurfaces $S_{2n-1}$ and $\overline{S}_{2n-1}$ are reducible for
$k = -2n + 1, \, \dots, 2n - 2$. The degree in $k$ of this map and its inverse is $(2n)^2$.
\end{conjecture}
Here the conjecture for the inverse maps:
\begin{conjecture} \label{ConjTwoSidedInv}
The determinant of the Jacobian $\mathrm{det} \,\mathrm{jac}$ of $g^{*-1} (k)$
\begin{equation}
  \mathrm{det}\,\mathrm{jac} \, (g^{*-1} (k)) \sim S_{2n-1} \; P^{n(2n-1)} \; \overline{P}^{n(2n-1)}
  \quad \text{of degree } (2n)^2 - 1
\label{fForm} \end{equation}
is the union of \textbf{three components}: a hypersurface $S_{2n-1}$ of degree $2n - 1$ and
multiplicity $1$ and two hyperplanes $P$ and $\overline{P}$, each of multiplicity $n (2n - 1)$.
The hyperplanes $P$ and $\overline{P}$ are permuted by the left-right symmetry (\ref{LR_Sym_j}).
\end{conjecture}
For $n = 1, 2$ the 2 conjectures above are true, see appendices \ref{SectBir2TwoSided} and \ref{SectBir4TwoSided}.
About the computation of these maps see appendix \ref{SectCAS}. \\
Further evidence for this conjecture based on symbolic calculations with a CAS for the cases
$n = 3, \, \dots$ i.e. maps in $Bir \left( \mathbb{P}_\mathbb{R}^6 \right), \, \dots$ is not available until
now due to limitations of my computing resources.

\section{Conclusion}

It is shown by symbolic computations, that there exists strong evidence for the conjectures
 \ref{ConjOneSided} and  \ref{ConjTwoSided} to be true. It would be interesting to get a proof of
the conjectures and so the existence of these two families, preferred without the use of an CAS.

\newpage

\noindent \textbf{\large Appendices}

\appendix

\section{Examples for the maps in the first, one-sided family \label{SectFirstOneSided}} 

\subsection{$n = 1$, the map $f^* (k)$ in the Cremona group $Bir \left( \mathbb{P}_\mathbb{R}^{1} \right)$ \label{SectBir1}} 

$$ F (k) = 1 + k (k + 1) / 2 \, l_0 $$
\begin{align}
f (k) : \enspace  & \mathbb{R}^1 & \rightarrow & \enspace \mathbb{RP}^{1} \nonumber \\
& \textbf{l} = (l_0) & \mapsto & \enspace ( \; F (k) : F (k + 1) \; )
 \label{MapTrans_c0}
\end{align}

\begin{proposition}
The map $f^* (k)$ is linear and therefore birational for $k \in \mathbb{R} \setminus \, \{-1 \}$.
The determinant of the Jacobian $\mathrm{det}\,\mathrm{jac}$ does not depend on \textbf{l},
it has degree~$0$.
\end{proposition}

\subsection{$n = 2$, the map $f^* (k)$ in the Cremona group $Bir \left( \mathbb{P}_\mathbb{R}^{2} \right)$ \label{SectBir2}} 

\begin{align} \begin{split}
E (k) & = 1 + (k + 1) (k + 2) (l_0 + 3 k (k + 3) (2 - (k - 1) (k + 1) (k + 2) (k + 4) \, l_1 )  \, l_1)  \\
F (k) & = 1 + k (k + 2) ((6 k^2 + 2 k - 1) + l_0) \, l_1 - 3 (k - 1) k (k + 1)^2 (k + 2) (k + 3) \, l_1^2)
\nonumber \end{split} \end{align}
\begin{align}
f (k) : \enspace  & \mathbb{R}^2 & \rightarrow & \enspace \mathbb{RP}^{2} \nonumber \\
& \textbf{l} = (l_0, l_1) & \mapsto & \enspace ( \; F (k) : E (k) : F (k + 1) \; )
 \label{MapTrans_c1}
\end{align}

\begin{proposition}
The map $f^* (k)$ is birational for $k \in \mathbb{R} \setminus \, \{-2, -3/2, -1 \}$.
The determinant of the Jacobian $\mathrm{det}\,\mathrm{jac}$ is the union of $2$ lines $L_1, L_2$
with multiplicity $1$ and $2$. Homogenized with the variable $h$, these 2 lines are:
$$L_1 : \enspace h = 0 \qquad L_2 : \enspace 3 \, k \, (k + 1) (k + 2) (k + 3) \, l_1 - h = 0$$
Let $\rho$ be the involution in example 4.1.1. in \cite{Des2012} $\rho: (l_0 : l_1 : h)
 \dashrightarrow (l_0 l_1 : h^2 : l_1 h)$, then $f_k^* = \alpha_k \circ \rho \circ \beta_k$ with
 $\alpha_k, \beta_k$ linear maps. \\
For $k = -3, \, \dots, 0$ the lines $L_1, L_2$ represent the same line, so we get a triple line.
\end{proposition}

\subsection{$n = 3$, the map $f^* (k)$ in the Cremona group $Bir \left( \mathbb{P}_\mathbb{R}^{3} \right)$ \label{SectBir3}} 

\begin{align} \begin{split}
E (k) & = 1 + (k + 3) (3 k + 2) / 6 \, l_0 + k (k + 3) (3 k^2 + 13 k + 9) / 12 \, l_1 \, \dots \\
F (k) & = 1 + k (k + 3) / 2 \, l_0 + k (k + 3) (k^2 + 3 k - 1) / 4 \, l_1 \, \dots 
\nonumber \end{split} \end{align}
$E (k)$ and $F (k)$ have $9$ monoms, the degree in $k$ of these functions is $18$.
\begin{align}
f (k) : \enspace  & \mathbb{R}^2 & \rightarrow & \enspace \mathbb{RP}^{2} \nonumber \\
& \textbf{l} = (l_0, l_1) & \mapsto & \enspace ( \; F (k) : E (k) : E (-k - 4) : F (-k - 4) \; )
 \label{MapTrans_c2}
\end{align}

\begin{proposition}
The map $f^* (k)$ is birational for $k \in \mathbb{R} \setminus \, \{ -3, -5/2, -2, -3/2, -1 \}$.
The determinant of the Jacobian $\mathrm{det}\,\mathrm{jac}$ is the union of a hyperplane $P$
with multiplicity $2$ and a hypersurface $S_2$ of degree $2$ and  multiplicity $3$.
Homogenized with the coordinate $h$, the hyperplane is $h = 0$.
The quadratic hypersurface $S_2$ is:
$$S_2 : \enspace (k - 1) k^2 (k + 1)^2 (k + 2)^2 (k + 3)^2 (k + 4)^2 (k + 5) \, l_2^2
 + h \, (\dots \hbox{ a linear form } \dots) = 0$$
\end{proposition}
From this form of the quadratic hypersurface we see that the quadratic factorizes, if the first
term above is $0$. Then the quadratic has the factor $h$ for $k \in \{ -5 \enspace (1) \enspace +1 \}$.
$\mathrm{det}\,\mathrm{jac}$ in these cases is the union of:
\begin{enumerate}[label=(\roman*)]
   \item the hyperplane $h = 0$ with multiplicity $5$ and a hyperplane $5 \, l_1 - 15 \, 1_2 - h = 0$
         with multiplicity $3$ for $k = -5, +1$.  
   \item the hyperplane $h = 0$ with multiplicity $8$ for $k = -4, \, \dots, 0$.  
\end{enumerate}

\subsection{$n = 4$, the map $f^* (k)$ in the Cremona group $Bir \left( \mathbb{P}_\mathbb{R}^{4} \right)$ \label{SectBir4}} 

\begin{align} \begin{split}
D (k) & = 1 + (2 k^2 + 10 k + 7) / 4 \, l_0 + (3 k^4 + 30 k^3 + 92 k^2 + 85 k + 21) / 12 \, l_1 \, \dots \\
E (k) & = 1 + (k + 4) (2 k + 1) / 4 \, l_0 + k (k + 4) (k^2 + 5 k + 3) / 4 \, l_1 \, \dots \\
F (k) & = 1 + k (k + 4) / 2 \, l_0 + k (k + 4) (k^2 + 4 k - 1) / 4 \, l_1 \, \dots 
\nonumber \end{split} \end{align}
$D (k), E (k)$ and $F (k)$ have $21$ monoms, the degree in $k$ of these functions is $32$.
\begin{align}
f (k) : \enspace  & \mathbb{R}^2 & \rightarrow & \enspace \mathbb{RP}^{2} \nonumber \\
& \textbf{l} = (l_0, l_1) & \mapsto & \enspace ( \; F (k) : E (k) : D (k) : E (-k - 5) : F (-k - 5) \; )
 \label{MapTrans_c3}
\end{align}

\begin{proposition}
The map $f^* (k)$ is birational for $k \in \mathbb{R} \setminus \, \{ -4, -7/2, -3, -5/2,$ $-2, -3/2, -1 \}$.
The determinant of the Jacobian $\mathrm{det}\,\mathrm{jac}$ is the union of a hyperplane $P$
with multiplicity $3$ and a hypersurface $S_3$ of degree $3$ and  multiplicity $4$.
Homogenized with the coordinate $h$, the hyperplane is $h = 0$.
The cubic hypersurface $S_3$ is:
\begin{align} \begin{split}
S_3 : \enspace & (k - 2) (k - 1)^2 k^3 (k + 1)^3 (k + 2)^3 (k + 3)^3 (k + 4)^3 (k + 5)^3
 (k + 6)^2 (k + 7) \, l_3^3 \\
 & + h \, (\dots \hbox{ a quadratic } \dots) = 0 
\nonumber \end{split} \end{align}
\end{proposition}
From this form of the cubic hypersurface we see that the cubic is reducible, if the first
term above is $0$. Then the cubic has the factor $h$ for $k \in \{ -7 \enspace (1) \enspace +2 \}$.
$\mathrm{det}\,\mathrm{jac}$ in these cases is the union of:
\begin{enumerate}[label=(\roman*)]
   \item the hyperplane $h = 0$ with multiplicity $7$ and a quadratic hypersurface
         with multiplicity $4$ for $k = -7, +2$.  
   \item the hyperplane $h = 0$ with multiplicity $11$ and a hyperplane
         $15 \, l_1 + 60 \, 1_2 + 135 \, l_3 - 2 \, h = 0$ with multiplicity $4$ for $k = -6, +1$.  
   \item the hyperplane $h = 0$ with multiplicity $15$ for $k = -5, \, \dots, 0$.  
\end{enumerate}

\subsection{$n = 2$, the \enspace i n v e r s e \enspace  map $f^{*-1} (k)$ in the Cremona group $Bir \left( \mathbb{P}_\mathbb{R}^{2} \right)$ \label{SectBir2OneSidedInv}} 

This map shows more structure than its counterpart $f^* (k)$. There exists a linear factor
$\Gamma (k)$ in $\textbf{j}$, that appears with increasing powers in the components. This
$\Gamma (k)$ in $\textbf{j}$ defines directly a line in the components of $\mathrm{det}\,\mathrm{jac}$. 
\begin{align} \begin{split}
  & P_1 (k), P_2 (k) \qquad \text{homogenous polynomials in } \textbf{j} \text{ of degree } 1, 2 \\
  & \Gamma (k) =      (k + 1) (k + 3) (k + 4) \, j_0 - k (k + 3) (2k + 3) \, j_1
                    + (k - 1) k (k + 2) \, j_2
\nonumber \end{split} \end{align}
\begin{align}
f^{*-1} (k) : \enspace  & \mathbb{RP}^2 & \rightarrow & \enspace \mathbb{RP}^2 \nonumber \\
& \textbf{j} = (j_0 : j_1 : j_2) & \mapsto & \enspace \textbf{l}, h = (l_0 : l_1 : \; h) = \label{MapTransInv_c1} \\
&    &    & \enspace ( \; P_2 (k)  \, \Gamma^0 (k) : P_1 (k) \, \Gamma^1 (k) : \nonumber \\
&    &    & \enspace \enspace  3 \, (k + 1) (k + 2) \, \Gamma^2 (k) \; )
\nonumber 
\end{align}

\begin{proposition}
The determinant of the Jacobian $\mathrm{det}\,\mathrm{jac}$ is the union of $2$ lines $L_1, L_2$
with multiplicity $2$ and $1$, these 2 lines are:
\begin{align} \begin{split}
  & L_1 : \Gamma (k) = 0 \quad L_2 : (k + 1) (k + 3) \, j_0 - k (k + 2) \, j_2 = 0
\nonumber \end{split} \end{align}
For $k = -3, -3/2, 0$ the lines $L_1, L_2$ represent the same line, so we get a triple line.
\end{proposition}

\subsection{$n = 3$, the \enspace i n v e r s e \enspace  map $f^{*-1} (k)$ in the Cremona group $Bir \left( \mathbb{P}_\mathbb{R}^{3} \right)$ \label{SectBir3OneSidedInv}} 

This map shows more structure than its counterpart $f^* (k)$. There exists a linear factor
$\Gamma (k)$ in $\textbf{j}$, that appears with increasing powers in the components.  This
$\Gamma (k)$ in $\textbf{j}$ defines directly the hyperplane in the components of $\mathrm{det}\,\mathrm{jac}$. 
\begin{align} \begin{split}
  & P_1 (k), P_2 (k), P_3 (k) \qquad \text{homogenous polynomials in } \textbf{j} \text{ of degree } 1, 2, 3 \\
  & \Gamma (k) =      (k + 1) (k + 4) (k + 5) (k + 6) (2k + 3) \, j_0 
                    - 3 k (k + 1) (k + 4) (k + 5) (2k + 5) \, j_1       \\
  & \qquad \quad    + 3 (k - 1) k (k + 3) (k + 4) (2k + 3) \, j_2
                    - (k - 2) (k - 1) k (k + 3) (2k + 5) \, j_3
\nonumber \end{split} \end{align}
\begin{align}
f^{*-1} (k) : \enspace  & \mathbb{RP}^3 & \rightarrow & \enspace \mathbb{RP}^3 \nonumber \\
& \textbf{j} = (j_0 : j_1 : j_2 : j_3) & \mapsto & \enspace \textbf{l}, h = (l_0 : l_1 : l_2 : \; h) = \label{MapTransInv_c2} \\
&    &    & \enspace ( \; P_3 (k)  \, \Gamma^0 (k) : P_2 (k) \, \Gamma^1 (k) : \nonumber \\
&    &    & \enspace \enspace  P_1 (k) \, \Gamma^2 (k) : (k + 1) (k + 2) (k + 3) \, \Gamma^3 (k) \; )
\nonumber 
\end{align}
\begin{proposition}
The determinant of the Jacobian $\mathrm{det}\,\mathrm{jac}$ of $f^{*-1} (k)$ is the union of the
hyperplane $P : \Gamma (k) = 0$ with multiplicity $6$ and a hypersurface $S_2$ of degree $2$ and  multiplicity $1$.  
\end{proposition}

\subsection{A reflection symmetry in the parameter $k$ of the map $f^* (k)$ \label{SectReflSymOne}} 

Let $\rho$ be the following linear, involutory Cremona map. This map reverses the components, it can
be viewed as a reflection at the middle of the components.
\begin{align}
\rho : \enspace  & \mathbb{RP}^{n} & \rightarrow & \enspace \mathbb{RP}^{n} \nonumber \\
& ( j_0 : j_1 : \, \dots : j_{n-1} : j_{n} ) & \mapsto &
 \enspace ( j_{n} : j_{n - 1} : \, \dots : j_1 : j_0 ) \nonumber
\end{align}
The following linear, involutory map represents a reflection at the point \\ $k_R = - (n + 1) / 2$:
\begin{align}
 R : \enspace  & \mathbb{R} & \rightarrow & \enspace \mathbb{R} \nonumber \\
& k & \mapsto & - n - 1 - k \label{MapReflkOne}
\end{align}
\begin{conjecture} \label{ConjOneSidedRefl}
$f^* (R (k)) = \rho \circ f^* (k)$
\end{conjecture}
This conjecture is true for $n = 1, 2, 3$ and $4$.

\section{Examples for the maps in the second, two-sided family \label{SectSecondTwoSided}} 

\subsection{$n = 1$, the map $g^* (k)$ in the Cremona group $Bir \left( \mathbb{P}_\mathbb{R}^{2} \right)$ \label{SectBir2TwoSided}} 

$$ H (k) = 1 + k^2 / 2 \, (l_0 + r_0 + (k^2 - 1) / 2 \, l_0 r_0) $$
\begin{align}
g (k) : \enspace  & \mathbb{R}^2 & \rightarrow & \enspace \mathbb{RP}^{2} \nonumber \\
& \textbf{l}, \textbf{r} = (l_0, \; r_0) & \mapsto & \enspace ( \; H (k) : (l_0 - r_0) : H (k + 1) \; )
\end{align}

\begin{proposition}
The map $g^* (k)$ is birational for $k \in \mathbb{R} \setminus \, \{ -1/2 \}$.
The determinant of the Jacobian $\mathrm{det}\,\mathrm{jac}$ is the union of $3$ lines
$L_1, L_2, \overline{L}_2$ each with multiplicity $1$. $L_2, \overline{L}_2$ are permuted by
the left-right symmetry (\ref{LR_Sym}). Homogenized with the coordinate $h$, these 3 lines are:
$$L_1 : \enspace h = 0 \qquad L_2 : \enspace k \, (k + 1) \, l_0 + 2 h = 0
 \qquad \overline{L}_2 : \enspace k \, (k + 1) \, r_0 + 2 h = 0$$
Let $\sigma$ be the standard involution \cite{Des2012} $\sigma: (l_0 : r_0 : h)
 \dashrightarrow (r_0 h : l_0 h : l_0 r_0)$, then $g_k^* = \alpha_k \circ \sigma \circ \beta_k$ with
 $\alpha_k, \beta_k$ linear maps. See \cite{BlHe2015}. \\
For $k = -1, 0$ the lines $L_1, L_2, \overline{L}_2$ represent the same line, so we get a triple line.
\end{proposition}

\subsection{$n = 2$, the map $g^* (k)$ in the Cremona group $Bir \left( \mathbb{P}_\mathbb{R}^{4} \right)$ \label{SectBir4TwoSided}} 

Though this map is quartic, depends on 1 real parameter and is in $Bir \left( \mathbb{P}_\mathbb{R}^{4} \right)$
it can be written in a not too lengthy form because of the left-right symmetry (\ref{LR_Sym}).
$$ H_0 (k, \textbf{d}) = 1 + (k - 1) k (d_0 + (k - 2) (k + 1) (6 - 3 (k - 3) (k - 1) k (k + 2) \, d_1) \, d_1) $$
\begin{align} \begin{split}
H (k) = \, & (  H_0 (k, \textbf{l}) H_0 (k + 1, \textbf{r})
              + H_0 (k, \textbf{r}) H_0 (k + 1, \textbf{l})) / 2 \\
           & - 36 \, (k^2 - 1) k^2 (l_1 - r_1)^2
\nonumber \end{split} \end{align}
$$ J_0 (k, \textbf{d}) = 1 + (k^2 + k + 3) \, d_0 + 6 \, (k^4 + 2 k^3 + k^2 + 6) \, d_1
                         - 3 \, (k^2 - 9) (k^2 - 4) (k^2 - 1) k (k + 4) \, d_1^2  $$
$$ J_1 (k, \textbf{d}) = 1 + k (k + 1) (d_0 + 3 \, (k - 1) (k + 2) (2 - (k - 2) k (k + 1) (k + 3) \, d_1) \, d_1) $$
\begin{align} \begin{split}
J (k) = \, & (J_0 (k, \textbf{l}) J_1 (k, \textbf{r}) + J_0 (k, \textbf{r}) J_1 (k, \textbf{l})) / 2 \\
           & + 108 \, (k^2 - 1) k (k + 2) (l_1 - r_1)^2
\nonumber \end{split} \end{align}
\begin{align} \begin{split}
K (k) = (l_0 - r_0) K_0 (k) K_0 (-k) \quad K_0 (k) = 1 - 3/2 \, (k^2 - 1) k (k + 2) \, (l_1 + r_1)
\nonumber \end{split} \end{align}
\begin{align} \begin{split}
L (k) = 3 \, (l_1 - r_1) \, n^2 & ( \quad 16 + 4 \, (k^2 - 1) (l_0 + r_0 - 4 \, (k^2 - 1) (l_1 + r_1)) \\
                                     & \; - 3 \, (k^2 - 4) (k^2 - 1)^2 \, (3 \, l_0 r_1 + 3 \, l_1 r_0 +
                                                                                l_0 l_1 + r_0 r_1 \\
                                     & \; \qquad \qquad \qquad \qquad \qquad + 16 \, (k^2 - 6) \, l_1 r_1 )  \\
                                     & ) \, / \, 4
\nonumber \end{split} \end{align}
\begin{align}
g (k) : \enspace  & \mathbb{R}^4 & \rightarrow & \enspace \mathbb{RP}^{4} \nonumber \\
& \textbf{l}, \textbf{r} = (l_0, l_1, \; r_0, r_1) & \mapsto & \enspace ( \; H (k) : (K (k) + L (k)) : J (k) :  \\
&                     &                 & \qquad \qquad \; (K (k + 1) + L (k + 1)) : H (k + 1) \; ) \nonumber 
\end{align}
$H (k)$ and $J (k)$ have $16$ monoms, the degree in $k$ of these functions is $16$. \\
$K (k) + L (k)$ has $12$ monoms, the degree in $k$ of this function is $10$.

\begin{proposition}
The map $g^* (k)$ is birational for $k \in \mathbb{R} \setminus \, \{ -3/2, -1/2, +1/2 \}$.
The determinant of the Jacobian $\mathrm{det}\,\mathrm{jac}$ is the union of one hyperplane $P$
with multiplicity $3$ and two hypersurfaces $S_3$ and $\overline{S}_3$ of degree $3$, each of
multiplicity $2$. Homogenized with the coordinate $h$, the hyperplane is $P : h = 0$. The cubics
$S_3, \overline{S}_3$ are permuted by the left-right symmetry (\ref{LR_Sym}).
One of the 2 cubic hypersurfaces is:
$$S_3 : \enspace (k - 2) (k - 1)^2 k^3 (k + 1)^3 (k + 2)^2 (k + 3) \, l_1^2 r_1
 + h \, (\dots \hbox{ a quadratic } \dots) = 0$$
\end{proposition}
From this form of the cubic hypersurface we see, that the cubic factorizes (if the first term above is $0$,
then the cubic has the factor $h$) for $k \in \{ -3, -2, -1, +0, +1, +2 \}$. \\
$\mathrm{det}\,\mathrm{jac}$ in these cases is the union of:
\begin{enumerate}[label=(\roman*)]
   \item the hyperplane $h = 0$ with multiplicity $7$ and two quadratic hypersurfaces each with
         multiplicity $2$ for $k = -3, +2$. 
   \item the hyperplane $h = 0$ with multiplicity $11$ and two hyperplanes $2 \, l_0 - 12 \, r_1
         + 12 \, l_1 + h = 0$ and $2 \, r_0 - 12 \, l_1 + 12 \, r_1 + h = 0$ each with
         multiplicity $2$ for $k = -2, +1$.  
   \item the hyperplane $h = 0$ with multiplicity $15$ for $k = -1, +0$.  
\end{enumerate}

\subsection{$n = 1$, the \enspace i n v e r s e \enspace  map $g^{*-1} (k)$ in the Cremona group $Bir \left( \mathbb{P}_\mathbb{R}^{2} \right)$ \label{SectBir2TwoSidedInv}} 

This map shows more structure than its counterpart $g^* (k)$. There exists a linear factor
$\Gamma (k)$ in $\textbf{j}$, this factor or its conjugate by (\ref{LR_Sym_j}) appear in all components.
This $\Gamma (k)$ in $\textbf{j}$ defines directly the 2 conjugated lines in the components of $\mathrm{det}\,\mathrm{jac}$.
\begin{align} \begin{split}
  & H (k)      =  2 \, j_0 - 2 (2k + 1) \, j_1 - 2 \, j_2 \\
  & \Gamma (k) =  2 (k + 1) (k + 2) \, j_0 + k (k + 1) (2k + 1) \, j_1 - 2 (k - 1) k \, j_2 
\nonumber \end{split} \end{align}
\begin{align}
g^{*-1} (k) : \enspace  & \mathbb{RP}^2 & \rightarrow & \enspace \mathbb{RP}^2 \nonumber \\
& \textbf{j} = (j_0 : j_1 : j_2 ) & \mapsto & \enspace \textbf{l}, \textbf{r}, h = (l_0 : \; r_0 : \; h) = \label{MapTransInv_c0} \\
&  &  & \enspace    ( \; 2 \, H (k) \, \Gamma (k) : 2 \,, \overline{H} (k) \, \overline{\Gamma} (k) :
                         \Gamma (k) \, \overline{\Gamma} (k)\; )
\nonumber 
\end{align}
\begin{proposition}
The determinant of the Jacobian $\mathrm{det}\,\mathrm{jac}$ of $g^{*-1} (k)$ is the union of
$3$ lines $L_1, L_2, \overline{L}_2$ each with multiplicity $1$. $L_2, \overline{L}_2$ are permuted
by the left-right symmetry (\ref{LR_Sym_j}). These 3 lines are:
\begin{align} \begin{split}
  L_1 : \enspace (k + 1) \, j_0 + k \, j_2 = 0 \qquad L_2 : \Gamma (k) = 0 \quad \overline{L}_2 :
  \overline{\Gamma} (k) = 0 
\nonumber \end{split} \end{align}
For $k = -1, -1/2, 0$ the lines $L_1, L_2, \overline{L}_2$ represent the same line, so we get a triple line.
\end{proposition}

\subsection{$n = 2$, the \enspace i n v e r s e \enspace  map $g^{*-1} (k)$ in the Cremona group $Bir \left( \mathbb{P}_\mathbb{R}^{4} \right)$ \label{SectBir4TwoSidedInv}} 

There exists a linear factor $\Gamma (k)$ in $\textbf{j}$, this factor or its conjugate by (\ref{LR_Sym_j})
appear with powers $1, 2$ in all components. This $\Gamma (k)$ in $\textbf{j}$ defines directly the 2 conjugated hyperplanes in the components of $\mathrm{det}\,\mathrm{jac}$.

\begin{align} \begin{split}
  & P_1 (k), P_2 (k) \qquad \text{homogenous polynomials in } \textbf{j}, \; 
    \text{of degree } 1, 2 \\
  & \Gamma (k) =      (k + 1) (k + 2) (k + 3) (k + 4) (2k - 1) \, j_0 
                    + k (k + 1) (k + 2) (k + 3) (2k - 1) (2k + 3) \, j_1       \\
  & \qquad \quad    - 2 (k - 1) k (k + 1) (k + 2) (2k + 1) \, j_2
                    - (k - 2) (k - 1) k (k + 1) (2k - 1) (2k + 3) \, j_3        \\
  & \qquad \quad    + (k - 3) (k - 2) (k - 1) k (2k + 3) \, j_4
\nonumber \end{split} \end{align}
\begin{align}
g^{*-1} (k) : \enspace  & \mathbb{RP}^4 & \rightarrow & \enspace \mathbb{RP}^4 \nonumber \\
& \textbf{j} = (j_0 : j_1 : j_2 : j_3 : j_4 ) & \mapsto & \enspace \textbf{l}, \textbf{r}, h = (l_0 : l_1 : \; r_0 : r_1 : \; h) = \label{MapTransInv_c1} \\
&  &  & \enspace    ( \;  P_2 (k) \, \Gamma^2 (k) : P_1 (k) \, \Gamma (k) \, \overline{\Gamma}^2 (k) :
\nonumber \\
&  &  & \enspace \enspace \overline{P}_2 (k) \, \overline{\Gamma}^2 (k) : \overline{P}_1 (k) \, \overline{\Gamma} (k) \, \Gamma^2 (k) :
                         9 \, \Gamma^2 (k) \, \overline{\Gamma}^2 (k)\; )
\nonumber 
\end{align}
\begin{proposition}
The determinant of the Jacobian $\mathrm{det}\,\mathrm{jac}$ of $g^{*-1} (k)$ is the union of a
hypersurface $S_{3}$ of degree $3$ and multiplicity $1$ and two hyperplanes $P : \Gamma (k) = 0$
and $\overline{P} : \overline{\Gamma} (k) = 0$, each of multiplicity $6$. The hyperplanes $P$ and
$\overline{P}$ are permuted by the left-right symmetry (\ref{LR_Sym_j}). 
\end{proposition}
For certain $k$ the cubic $S_{3}$ is reducible, $\mathrm{det}\,\mathrm{jac}$ in these cases is the union of:
\begin{enumerate}[label=(\roman*)]
   \item a hyperplane $j_0 - j_4 = 0$ with multiplicity $3$ and two conjugated hyperplanes, each
         with multiplicity $6$ for $k = -1/2$.  
   \item a hyperplane with multiplicity $13$ and a second hyperplane with multiplicity $2$
         for $k = -3/2, +1/2$.  
   \item a hyperplane with multiplicity $15$ for $k = -1, 0$.  
\end{enumerate}

\subsection{A reflection symmetry in the parameter $k$ of the map $g^* (k)$ \label{SectReflSymTwo}} 

Let $\rho$ be the following linear, involutory Cremona map. This map reverses the components, it can
be viewed as a reflection at the middle of the components.
\begin{align}
\rho : \enspace  & \mathbb{RP}^{2n} & \rightarrow & \enspace \mathbb{RP}^{2n} \nonumber \\
& ( j_0 : j_1 : \, \dots : j_{2n-1} : j_{2n} ) & \mapsto &
 \enspace ( j_{2n} : j_{2n-1} : \, \dots : j_1 : j_0 ) \nonumber
\end{align}
The following linear, involutory map represents a reflection at the point $k_R = - 1 / 2$:
\begin{align}
 R : \enspace  & \mathbb{R} & \rightarrow & \enspace \mathbb{R} \nonumber \\
& k & \mapsto & - 1 - k \label{MapReflkTwo}
\end{align}
\begin{conjecture} \label{ConjTwoSidedRefl}
$f^* (R (k)) = \rho \circ f^* (k)$
\end{conjecture}
This conjecture is true for $n = 1$ and $2$.

\section{Computations done with a computer algebra system \label{SectCAS}}

The maps in the appendices \ref{SectFirstOneSided} \dots \ref{SectSecondTwoSided} were computed with
the aid of the CAS Maxima. Here I give a very rough sketch how the polynomials and maps are computed.
The sketch of the algorithm is done here for the $Q_i (\textbf{l}, x)$ and $f_i$ in the first family,
but the computation of the $M_i (\textbf{l}, \textbf{r}, x)$ and $g_i$ in the second family can be done
in a similar manner: \\

\begin{enumerate}[label=\arabic*.]
   \item step : \underline{calculate all $Q_i (\textbf{l}, x)$ up to a certain index} \\
         set $Q_0 (\textbf{l}, x) = 1$, calculate the further $Q$'s by recursion. Assuming all
         $Q_i (\textbf{l}, x)$ are calculated up to the index k, the recursion step to calculate
         $Q_{k+1} (\textbf{l}, x)$ is: 
  \begin{enumerate}[label=\alph*)]
     \item make the "Ansatz" (\ref{QasSumJ}) for $Q_{k+1} (\textbf{l}, x)$ which expresses the $Q$
           as a sum of symmetric orthogonal Jacobi polynomials. The Jacobi polynomials $J_i (x)$ with
           negative index $i$ are set to $0$. The unknowns $j_{{k+1},0} (\textbf{l})$
           \, \dots \, $j_{{k+1},n} (\textbf{l})$ have to be determined 
     \item $Q_{k+1} (\textbf{l}, x)$ has to be orthogonal (\ref{OrthQ}) to the already calculated
           polynomials with lower indices, we get $k + 1$ linear equations in the $n + 1$
           unknowns $j_{{k+1},i} (\textbf{l})$ 
     \item solve this system of homogeneous linear equations 
     \item the just $n + 1$ calculated $j_{{k+1},i} (\textbf{l})$ have denominators in $\textbf{l}$,
           multiply it by the least common multiple of all denominators. Now the $j_{{k+1},i} (\textbf{l})$
           are polynomials in $\textbf{l}$
     \item the new $Q_{k+1} (\textbf{l}, x)$ is now determined by (\ref{QasSumJ}) 
  \end{enumerate}
   \item step : \underline{get the $j_{k,i} (\textbf{l})$ as polynomials in $k$} \\
         in the previous step we calculated the $j_{k,i} (\textbf{l})$ for integer $k$. The coefficients
         of the monoms in $\textbf{l}$ are rational numbers. The get 1-parameter maps with
         $k \in \mathbb{R}$ instead of only $k \in \mathbb{N}_+$ we need these coefficients as functions
         of $k$. Therefore the coefficient of each monom in $\textbf{l}$ has to be expressed as a
         polynomial in $k$ obtained by interpolating for various $k$ with $n \le k \, \dots$.
         How many of these $k$ we need depends on the degree in $k$ of the interpolation. This degree is
         conjectured to be $2n^2$, see conjecture \ref{ConjOneSided}.  \\
         Because the system of linear equations in step 1. is homogeneous we have to interpolate
         \emph{ratios} of the coefficients, take as denominator in this ratio the coefficient of the
         monom $h^n$, $h$ is the homogenizing variable
   \item step : \underline{beautify the maps} \\
         factorize the coefficients of the monoms in $\textbf{l}$, the polynomials in $k$, use brackets
         and common subfactors of the monoms, use symmetries \dots  
\end{enumerate}

\section{The connection to quadrature rules for splines \label{SectNumAna}}

Originally the maps and polynomials in this article were found in the context of quadrature rules
for polynomial splines, a topic in Numerical Analysis. \\

Here I give a very short and informal description of quadrature rules for splines and the connection
to the Cremona maps in this article: \\
Given a closed interval $[a, b]$, this is subdivided in non-overlapping subintervals,
the boundaries of these subintervals are called knots. \\
A spline $s (x)$ is a piecewise (on the subintervals) polynomial function, the spline is of continuity
class $c$ when on the inner knots the spline up to the $c^{th}$ derivation is continuous. \\
We want to determine nodes $x_i \in [a, b]$ and weights $w_i$ so that the weighted summation
$\sum_{}^{} w_i s (x_i)$ is equals to the integral $\int_{a}^{b} s (x) \, dx$ for splines up to a
certain degree. Rules are optimal or suboptimal if they use the minimum number of nodes. \\

The continuity class of the splines is $c = n - 1$. With the map $f_k$ another map called recursion map
from $\mathbb{R}^{c+1} \rightarrow \mathbb{R}^{c+1}$ is defined by:
\begin{equation} Rec_k  = C \circ \underbrace { f_k^{-1} \circ r \circ f_k }_{R_k} = C \circ R_k
\label{RecMap} \end{equation}
$C, r$ are monominal involutions, $R_k$ an involution, so $Rec_k$ is a product of 2 involutions.

This $Rec_k$ allows us to get 2 subsequent subintervals with the nodes defined by the roots of the polynomials $Q_k (\textbf{l}, x)$ and $Q_k (Rec_k (\textbf{l}), x)$, $\textbf{l}$ arbitrary. Another
possibility for 2 subsequent subintervals is $Q_k (\textbf{l}, x)$ and $M_{k+1} (Rec_k (\textbf{l}),  \textbf{r}, x)$, $\textbf{l}, \textbf{r}$ arbitrary. With the corresponding weights
the summation is exact for all splines with a support of 1 or 2 subintervals up to a certain degree.
The exactness for splines with a support of 2 subintervals is achieved by applying the recursion
map $Rec_k$ to the $\textbf{l}$ argument. \\

\bibliographystyle{amsplain}

\end{document}